\documentclass{amsart}
\usepackage{amssymb} 
\usepackage{amsmath} 
\usepackage{amsbsy}
\usepackage{hyperref}
\usepackage{verbatim}
\usepackage{tikz}
\usepackage{tabu}
\usepackage{longtable}
\usepackage{multirow}

\DeclareMathOperator{\ord}{\upsilon}

\begin {document}

\newtheorem{theorem}{Theorem}
\newtheorem{lemma}{Lemma}[section]

\title[]{On perfect powers that are sums of two Fibonacci numbers}

\author{Florian Luca}
\address{School of Mathematics, University of the Witwatersrand,
Private Bag 3, Wits 2050, Johannesburg, South Africa; 
Max Planck Institute for Mathematics,
Vivatsgasse 7, 53111 Bonn, Germany and
Department of Mathematics, Faculty of Sciences, University of Ostrava,
30 dubna 22, 701 03 Ostrava 1, Czech Republic}
\email{florian.luca@wits.ac.za}

\author{Vandita Patel}
\address{Department of Mathematics, University of Toronto, Bahen Centre, 40 St. George St., Room 6290, Toronto, Ontario, Canada, M5S 2E4}
\email{vandita\_patel@hotmail.co.uk}

\thanks{
The first-named author
is supported by  grants CPRR160325161141 and an A-rated researcher award both from the NRF of South Africa and by grant no. 17-02804S of the Czech Granting Agency. 
}

\date{\today}

\keywords{Exponential equation, Hemachandra numbers, Fibonacci numbers}
\subjclass[2010]{Primary 11D61, Secondary 11B39}

\begin {abstract}
We study the equation $F_n + F_m = y^p$, where $F_n$ and $F_m$ are respectively the $n$-th and $m$-th Fibonacci numbers
and $p \ge 2$. We find all solutions under the assumption $n \equiv m \pmod{2}$.
\end {abstract}
\maketitle

\section{Introduction}
Fibonacci numbers are prominent as well as being ancient. Their first
known occurrence dates back to around $700$AD, in the work of
Virah\={a}\.{n}ka. Virah\={a}\.{n}ka's original work has been lost, but
is nevertheless cited clearly in the work of Gop\={a}la (c. $1135$); below is a
translation of \cite[pg. 101]{HC};

\medskip

\begin{quote}
``For four, variations of meters of two [and] three being mixed, five happens. For five, variations of two earlier Ð three [and] four, being mixed, eight is obtained. In this way, for six, [variations] of four [and] of five being mixed, thirteen happens. And like that, variations of two earlier meters being mixed, seven morae [is] twenty-one. In this way, the process should be followed in all m\={a}tr\={a}--v\d{r}ttas."
\end{quote}

\medskip

The sequence is discussed rigorously in the work of Jain scholar Acharya Hemachandra (c. $1150$, living in what is known today as Gujarat) about $50$ years earlier than Fibonacci's {\it Liber Abaci} ($1202$). Hemachandra was in fact studying {\it Sanskrit prosody} (poetry meters and verse in Sanskrit) and not mathematics. Given a verse with an ending of $n$ beats to fill, where the choice of beats consists of length $1$ (called {\it short}) and length $2$ (called {\it long}), in how many ways can one finish the verse? The answer lies within the fundamental sequence, defined by the recurrence; 
\begin{equation}\label{eq:HC}\tag{$\diamondsuit$}
H_{n+2} = H_{n+1} + H_{n},\qquad  H_1 = 1, \quad  H_2 = 2, \quad n\geq 1,
\end{equation}
where Hemachandra makes the concise argument that any verse that is to be
filled with $n$ beats must end with a long or a short beat. Therefore, this
recurrence is enough to answer the question: given a verse with $n$ beats
remaining, one has $H_n$ ways of finishing the {\it prosody}, with $H_n$
satisfying \eqref{eq:HC}.

\medskip

Since the $12$th century, the Hemachandra/Fibonacci numbers have sat in the spotlight of modern number theory. They have been vastly studied; intrinsically for their beautiful identities but also for their numerous applications, for example, the golden ratio has a regular appearance in art, architecture and the natural world!

\medskip

Finding all perfect powers in the Fibonacci sequence was a fascinating long-standing conjecture. In 2006, this problem was completely solved by Y. Bugeaud, M. Mignotte and S. Siksek (see \cite{Ann}), who innovatively combined the  modular approach with classical linear forms in logarithms. In addition to this, Y. Bugeaud, F. Luca, M. Mignotte and S. Siksek also found all of the integer solutions to 
\begin{equation}\label{eq:Fib1}
F_n \pm 1 = y^p \quad p\geq 2,
\end{equation}
(see \cite{BLMS}). The authors found a clever factorisation which descended the problem to finding solutions of $F_n  = y^p$. 

\medskip

In this paper, we consider the natural generalisation, 
\begin{equation}\label{eq:FibSum}
F_n \pm F_m = y^p, \quad p\geq 2. 
\end{equation}
%
%

\begin{theorem}
\label{thm:1}
All solutions of the Diophantine equation~\eqref{eq:FibSum}
in integers $(n,m,y,p)$ with $n \equiv m \pmod{2}$ either have $\max\{|n|,|m|\}\le 36$, or $y=0$ and $|n|=|m|$.  
\end{theorem} 
Since $F_1=F_2=1$, it follows that every solution $(n,y,p)$ of equation~\eqref{eq:Fib1} can be thought of as a solution 
$(n,m,y,p)$ of equation~\eqref{eq:FibSum} with $m=1$, $2$ according to whether $n$ is odd or even.
Therefore, Theorem~\ref{thm:1} is a genuine generalisation of the main result from \cite{BLMS}. 

For a complete list of solutions to equation~\eqref{eq:FibSum} with $\max\{|n|,|m|\} \le 1000$ without the parity 
restriction, see Section~\ref{sec:conjecture}.

\section*{Acknowledgement}
The authors would like to thank the Max Planck Institute for Mathematics, for their generous hospitality and provisions of a fantastic collaborative working environment.

\section{Preliminaries}

Let $\left(F_n\right)_{n\ge 0}$ be the Hemachandra/Fibonacci sequence given by;
\[
F_{n+2} = F_{n+1} + F_{n},\qquad  F_0 = 0, \quad  F_1 = 1, \quad n\geq 0.
\] 
Recall that $(F_n)_{n\ge 0}$ can be extended to be defined on the negative indices by using the above recurrence and giving $n$ the values $n=-1,-2,\ldots$. Thus the formula $F_{-n}=(-1)^{n+1}F_n$ holds for all $n$. 

\medskip

Let $(L_n)_{n\ge 0}$ be the Lucas companion sequence of the Hemachandra/Fibonacci sequence given by;
\[
L_{n+2} = L_{n+1} + L_{n},\qquad  L_0 = 2, \quad  L_1 = 1, \quad n\geq 0.
\] 
Similarly, this can also be extended to negative indices $n$, and the formula $L_n=(-1)^n L_{-n}$ holds for all $n$. 

\medskip
The Binet formulas for $F_n$ and $L_n$ are;
\begin{equation}
\label{eq:binet}
F_n=\frac{1}{\sqrt{5}} (\alpha^n-\beta^n)\quad \text{ and }\quad L_n=\alpha^n+\beta^n\quad {\text{\rm for~all}}\quad n\in {\mathbb Z},
\end{equation}
where $(\alpha,\beta)=((1+{\sqrt{5}})/2, (1-{\sqrt{5}})/2)$. There are many formulas relating Hemachandra/Fibonacci numbers and Lucas numbers. Two of which are useful for us are;
\begin{equation}
\label{eq:X} 
F_{2n}=F_n L_n \quad \text{ and }\quad L_{3n}=L_n(L_n^2-3(-1)^n),
\end{equation} 
which hold for all $n$. They can be proved using Binet's formulae \eqref{eq:binet}.

The following result is well-known and can also be proved using Binet's formulae \eqref{eq:binet}.

\begin{lemma}
\label{lem:0}
Assume $n\equiv m\pmod 2$. Then 
\[
F_n+F_m=
\begin{cases} F_{(n+m)/2} L_{(n-m)/2} & \text{ if }\;  n\equiv m\pmod 4,\\
F_{(n-m)/2} L_{(n+m)/2} & \text{ if }\;  n\equiv m+2\pmod 4.
\end{cases}
\]
Similarly, 
\[
F_n-F_m=
\begin{cases} F_{(n-m)/2} L_{(n+m)/2} & \text{ if }\;  n\equiv m\pmod 4,\\
F_{(n+m)/2} L_{(n-m)/2} & \text{ if }\;  n\equiv m+2\pmod 4.
\end{cases}
\]
\end{lemma}

 The following result can be found in \cite{Md}.

\begin{lemma}
\label{lem:1}
Let $n=2^a n_1$ and $m=2^bm_1$ be positive integers with $n_1$ and $m_1$ odd integers and $a$ and $b$ nonnegative integers. Let $d=\gcd(n,m)$. Then
\begin{itemize}
\item[i)] $\gcd(F_n,F_m)=F_d$.
\item[ii)] $\gcd(L_n,L_m)=L_d$ if $a=b$ and it is $1$ or $2$ otherwise.
\item[iii)] $\gcd(F_n,L_m)=L_d$ if $a>b$ and it is $1$ or $2$ otherwise.   
\end{itemize}
\end{lemma}

The following results can be extracted from \cite{Crell}, \cite{JNTB} and \cite{Ann} and will be useful for us.

\begin{theorem}
\label{lem:2} 
If 
$$
F_n=2^s\cdot y^b
$$
for some integers $n\ge 1,~y\ge 1,~b\ge 2$ and $s \ge 0$ then $n\in \{1,2,3,6,12\}$. The solutions of the similar equation with $F_n$ replaced by $L_n$ have $n\in \{1,3,6\}$.  
\end{theorem}
\begin{theorem}
\label{lem:3} 
If 
$$
F_n=3^s\cdot y^b
$$
for some integers $n\ge 1,~y\ge 1,~b\ge 2$ and $s \ge 0$ then $n \in \{1,2,4,6,12\}$. The solutions of the similar equation with $F_n$ replaced by $L_n$ have $n \in \{1,2,3\}$.  
\end{theorem}

The following result is due to McDaniel and Ribenboim (see \cite{RM}).

\begin{theorem}
\label{thm:RM}
\begin{itemize}
\item[i)] Assume $u\mid v$ are positive integers such that $F_v/F_u=y^2$. Then, either $u=v$ or $(v,u)\in \{(12,1),~(12,2),(2,1),~(6,3)\}$.  
\item[ii)] Assume that $u\mid v$, $v/u$ is odd and $L_v/L_u=y^2$. Then, $u=v$ or $(v,u)=(3,1)$.
\end{itemize}
\end{theorem}

\section{Perfect Powers from Products of a Fibonacci and a Lucas Number}

\begin{theorem}\label{thm:FNxLM}
The only solutions to 
\[
F_N \cdot L_M=2^s\cdot y^p
\]
with $N$, $M$, $y$ positive integers, $s \ge 0$ and $p \ge 2$ satisfy 
\begin{gather*}
(N,M)=( 1 , 1 ), \quad
( 1 , 3 ), \quad
( 1 , 6 ), \quad
( 2 , 1 ), \quad
( 2 , 3 ), \quad
( 2 , 6 ), \quad
( 3 , 1 ), \quad
( 3 , 3 ), \\
( 3 , 6 ), \quad
( 4 , 2 ), \quad
( 4 , 6 ), \quad
( 6 , 1 ), \quad
( 6 , 3 ), \quad
( 6 , 6 ), \quad
( 12 , 1 ), \quad
( 12 , 2 ), \\
( 12 , 3 ), \quad
( 12 , 6 ), \quad
( 24 , 12 ). 
\end{gather*}
\end{theorem}
\begin{proof}
We shall in fact show that $N \le 24$ and $M \le 12$. The proof is then completed by a simple
program.
Write
\[
N=2^a N_1, \qquad M=2^b M_1,
\]
where $N_1$, $M_1$ are odd.
If $a \le b$, then by Lemma~\ref{lem:1}, we know $\gcd(F_N,L_M)=1$ or $2$, so
$F_N=2^u y_1^p$ and $L_M=2^v y_2^p$. By Theorem~\ref{lem:2}, we deduce that $N \le 12$ and $M \le 6$. 

Thus, we may assume that $a>b$. Let $r=a-b \ge 1$ and $d=\gcd(N,M)$. Therefore, $d=2^b \gcd(N_1,M_1)$. Write $N=2^r k d$
where $k$ is odd. Then we obtain;
\[
2^s y^p= F_N \cdot L_M=F_{2^r k d}\cdot L_M=F_{kd} \cdot L_{kd} \cdot L_{2kd} \cdots L_{2^{r-1}kd} \cdot L_M,
\]
by repeated application of \eqref{eq:X}. Note that 
\[
\ord_2(kd)=\ord_2(M), \qquad \text{$\ord_2(kd) \le \ord_2(2^i kd)$ \quad for $i \ge 0$}. 
\]
Thus, by Lemma~\ref{lem:1}, the greatest common divisor of $F_{kd}$ and $L_{kd} \cdot L_{2kd} \cdots L_{2^{r-1}kd} \cdot L_M$
is a power of $2$. Hence,
\[
F_{kd}=2^u y_1^p.
\]
By Theorem~\ref{lem:2}, we take note that $kd \in \{1,2,3,6,12\}$. Moreover,
\[
L_{kd} \cdot L_{2kd} \cdots L_{2^{r-1}kd} \cdot L_M=2^v y_2^p.
\]
Suppose $r \ge 2$. Then $\ord_2(2^{r-1}kd)> \ord_2(M)$. Once more, we use Lemma~\ref{lem:1} to see that  
the greatest common divisor of $L_{2^{r-1}kd}$ and $L_{kd} \cdot L_{2kd} \cdots L_{2^{r-2}kd} \cdot L_M$ is a power of $2$. Hence,
\[
L_{2^{r-1} kd}=2^w y_3^p.
\]
By Theorem~\ref{lem:2}, we conclude that $2^{r-1} kd \in \{1,3,6\}$. Therefore, $kd=3$ and $r=1$, contradicting our
assumption that $r \ge 2$. Hence, $r=1$ and $N \in \{2,4,6,12,24\}$. 

If $N=2$ or $6$ then $F_N=1$ or $8$ so $L_M=2^{\sigma} y^p$. Theorem~\ref{lem:2} allows us to readily conclude that $M \le 6$.
The cases $N=4$, $12$ and $24$ remain and require delicate treatment. 
\medskip
First, we deal with the cases $N=4$ and $N=12$. Since $F_4=3$ and $F_{12}=2^4 \times 3^2$
we have $L_M=2^\alpha 3^{\beta} y_0^p$ where $y_0$ is odd. 
If $\alpha=0$ then by Theorem~\ref{lem:3} we know that $M \le 3$ and so we may suppose that $s \ge 1$.
Thus, $2 \mid L_M$. Note that $2 \mid\mid M$ (as $r=1$ and $a=2$).
As $6 \mid L_M$, we have $3 \mid M$. Thus, we can write $M=2 \cdot 3^t \cdot \ell$, where $\ell$ is coprime to $6$.
Now, observe that
\[
L_{3\delta}=L_\delta (L_\delta^2-3)
\]
for $\delta=2 \cdot 3^i \cdot \ell$ with $i=t-1,t-2,\dots,0$. Hence, we obtain inductively $L_{2 \cdot 3^i \cdot \ell}=2^{s_i} 3^{w_i} y_i^p$.
For $i=0$, we have $3 \nmid 2 \ell$ and so $s_0=0$. Using Theorem~\ref{lem:3}, we infer that $\ell=1$. Hence, $M=2 \cdot 3^t$.
Note that $107 \mid\mid L_{18} \mid L_M$ and $107^2 \nmid L_{n}$ unless $(18 \times 107) \mid n$.  We conclude that $t \le 1$ and 
so $M \le 6$.

\medskip
Finally, let $N=24$, whence $F_{24}=2^5 \cdot 3^2 \cdot 7 \cdot 23$. Moreover, $M=2^2 \cdot M_1$ where $M_1$ is odd.
Thus, $d=2^2 \cdot \gcd(3,M_1)=4$ or $12$. However, 
$\gcd(F_N,L_M)=L_d$. 
We may rewrite the equation $F_N \cdot L_M=2^s y^p$
as
\[
L_d^2 \cdot \frac{F_N}{L_d} \cdot \frac{L_M}{L_d}=2^s y^p.
\]
If $d=4$ then $L_d=7$ and
we see that $23$ divides the left-hand side exactly once, giving a contradiction.
Thus, $d=12$ and so $L_d=2 \times 7 \times 23$. In this case the left-hand side is divisible by $3$ exactly twice
and therefore $p=2$. Thus, $L_M/L_{12}=2^\alpha y_1^2$. Since $M$ is an odd multiple of $12$, we can easily see that $L_M/L_{12}$
is odd and so $\alpha=0$. Hence, we can apply Theorem~\ref{thm:RM} to draw the inference that $M=12$.
\end{proof}

\section{Proof of Theorem~\ref{thm:1}}
If either $n=0$ or $m=0$, then the theorem follows from \cite{Ann}.
Via the identity $F_{-n}=(-1)^{n+1} F_n$, we can suppose that $n \ge m >0$. Note that changing
signs does not change parities, so we maintain the assumption $n \equiv m \pmod{2}$.


If $n=m$, then we need to solve $2F_n = y^p$, which is equivalent to solving $F_n = 2^{p-1}y_1^p$.
By Theorem~\ref{lem:2}, we have that $n \leq 12$.

Thus, we may suppose that $n>m > 0$.  By Lemma~\ref{lem:0}, there is some $\varepsilon = \pm 1$ such that letting
$N=(n+\varepsilon m)/2$ and $M=(n-\varepsilon m)/2$, we have $F_n \pm F_m=F_N \cdot L_M$. Observe that
$N$ and  $M$ are both positive. By Theorem~\ref{thm:FNxLM},
we know that $N \le 24$ and $M \le 12$. We finally conclude that  $n=N+M \le 36$ and $m=|N-M| <24$. This completes the proof.

\section{An Open Problem}\label{sec:conjecture}

It is still an open problem to find all solutions to equation~\eqref{eq:FibSum}
in the case $n \not\equiv m \pmod{2}$.
Under the condition, $n\not\equiv m \pmod{2}$, no factorisation is known for the left-hand side.
We searched for solutions with $0 \le m \le n \le 1000$ and found the following:
\begin{gather*}
F_0+F_0=0, \quad
F_1+F_0=1, \quad
F_2+F_0=1, \quad
F_3+F_3=2^2, \quad
F_4+F_1=2^2, \\
F_4+F_2=2^2, \quad
F_5+F_4=2^3, \quad
F_6+F_0=2^3,\quad
F_6+F_1=3^2,\quad
F_6+F_2=3^2,\\
F_6+F_6=2^4 = 4^2,\quad
F_7+F_4=2^4 = 4^2, \quad
F_9+F_3=6^2,\quad
F_{11}+F_{10}=12^2, \\
F_{12}+F_0=12^2,\quad
F_{16}+F_7=10^3,\quad
F_{17}+F_4=40^2,\quad
F_{36}+F_{12}= 3864^2.\\
\\
F_n-F_n=0, \quad F_1-F_0=1, \quad F_2-F_0=1, \quad F_2-F_1=0,\quad
F_3-F_1=1,\\
F_3- F_2=1,\quad
F_4- F_3=1,\quad
F_5- F_1=2^2,\quad
F_5- F_2=2^2,\quad
F_6- F_0=2^3,\\
F_7- F_5=2^3,\quad
F_8- F_5=2^4 = 4^2,\quad
F_8- F_7=2^3,\quad
F_9- F_3=2^5,\\
F_{11}- F_6=3^4 = 9^2,\quad
F_{12}- F_0=12^2,\quad
F_{13}- F_6=15^2,\quad
F_{13}- F_{11}=12^2,\\
F_{14}- F_{9}=7^3,\quad
F_{14}- F_{13}=12^2,\quad
F_{15}- F_{9}=24^2.
\end{gather*}
We conjecture that the above lists all of the solutions to equation~\eqref{eq:FibSum} with the restriction $n \ge m \ge 0$.

\end{document}